\documentclass[conf]{new-aiaa}
\usepackage[utf8]{inputenc}

\usepackage{graphicx}
\usepackage{amsmath}
\usepackage[version=4]{mhchem}
\usepackage{siunitx}
\usepackage{longtable,tabularx}
\usepackage{multirow}
\usepackage{wrapfig}
\usepackage{makeidx}
\usepackage{subfigure}
\usepackage{booktabs}
\title{Assessing the Performance of an Adaptive Multi-Fidelity Gaussian Process with Noisy Training Data: \\ A Statistical Analysis}

\author{S. Ficini\footnote{Ph.D. student, Department of Engineering, Roma Tre University, Via Vito Volterra 62, 00146 Rome. Email: simone.ficini@uniroma3.it} and U. Iemma\footnote{Professor, Department of Engineering, Roma Tre University, Via Vito Volterra 62, 00146 Rome.}}
\affil{Roma Tre University, Department of Engineering, Rome, Italy}
\author{R. Pellegrini\footnote{Postdoctoral Research Fellow, INM Rome, Via di Vallerano 139, 00128, AIAA Member.}, A. Serani\footnote{Research Scientist, INM Rome, Via di Vallerano 139, 00128.} and M. Diez\footnote{Senior Research Scientist, INM Rome, Via di Vallerano 139, 00128, AIAA Member. Email: matteo.diez@cnr.it}}
\affil{CNR-INM, National Research Council-Institute of Marine Engineering, Rome, Italy}

\begin{document}

\maketitle

\begin{abstract}
Despite the increased computational resources, the simulation-based design optimization (SBDO) procedure can be very expensive from a computational viewpoint, especially if high-fidelity solvers are required. Multi-fidelity metamodels have been successfully applied to reduce the computational cost of the SBDO process. In this context, the paper presents the performance assessment of an adaptive multi-fidelity metamodel based on a Gaussian process regression (MF-GPR) for noisy data. The MF-GPR is developed to: (i) manage an arbitrary number of fidelity levels, (ii) deal with objective function evaluations affected by noise, and (iii) improve its fitting accuracy by adaptive sampling. Multi-fidelity is achieved by bridging a low-fidelity metamodel with metamodels of the error between successive fidelity levels. The MF-GPR handles the numerical noise through regression. The adaptive sampling method is based on the maximum prediction uncertainty and includes rules to automatically select the fidelity to sample. The MF-GPR performance are assessed on a set of five analytical benchmark problems affected by noisy objective function evaluations. Since the noise introduces randomness in the evaluation of the objective function, a statistical analysis approach is adopted to assess the performance and the robustness of the MF-GPR. The paper discusses the efficiency and effectiveness of the MF-GPR in globally approximating the objective function and identifying the global minimum. One, two, and three fidelity levels are used. The results of the statistical analysis show that the use of three fidelity levels achieves a more accurate global representation of the noise-free objective function compared to the use of one or two fidelities.

\end{abstract}



\section{Introduction}

\lettrine{I}{n} the context of engineering design problems, where time and computational resources are usually limited, the simulation-based design (SBD) procedure has proven its ability to help designers in exploring design spaces and provides a large sets of design options. 
The SBD procedure explores design spaces and provides a large sets of design options by assessing design performance also for large sets of operating and environmental conditions \cite{serani2021hull}. The continuous development of high-performance computing systems has driven the SBD towards the automatic integration with optimization algorithms. The simulation-based design optimization (SBDO) procedure efficiently and effectively combines: (i) design modification and automatic meshing methods, (ii) numerical solvers, and (iii) optimization algorithms. The SBDO procedure has proven its ability, robustness and reliability to help designers in achieving global optimal design solutions. 

For the design of innovative configurations and off-design condition usually are required high-fidelity solvers to assess the accurate prediction of the desired output. Furthermore, to converge to the optimal solution the use of global optimization algorithms usually requires a large number of function evaluations. Despite the increased computational resources, the SBDO procedure with high-fidelity solvers can be very expensive from a computational viewpoint.

Furthermore, potential design improvement depends on the dimension and the extension of the design space. To achieve bigger improvement in finding the global optimum solution, high dimensionality and variability space need to be explored \cite{d2020design}. When the number of design variables increases the algorithms’ complexity and computational cost rapidly increase. Therefore, the global optimization are affected by the \textit{curse of dimensionality}, when the computational cost and the algorithms’ complexity increase with the problem dimension. A wide variety of techniques such as the design-space dimensionality reduction \cite{d2020design} can be used in order for the SBDO procedure to reduce the curse of dimensionality and be successful to efficiently describe the design variability with as few variables as possible.

Finally, global optimization algorithms ensure global exploration of the design space, but, to converge to the optimal solution, they usually require a large number of function evaluations. Therefore, the identification of the global design optimum is a challenging task. One methodology that has been developed to reduce the computational cost of the SBDO procedure, as here described, is metamodelling. Metamodels have been developed and successfully applied in several engineering fields \cite{viana2014-AIAA, serani2021hull} to reduce the computational cost of the optimization process. Using few simulations, a metamodel can provide an approximate and inexpensive to evaluate model of expensive simulations. Among others, Gaussian process (GP) \cite{coppede2019hydrodynamic}, dynamic Kriging \cite{zhao2010metamodeling}, and stochastic radial basis functions (SRBF) \cite{volpi2015development} have been successfully applied in the SBDO context. The research of accurate and efficient methods for metamodel-based analysis and optimization has moved from standard (or static) to function-adaptive approaches, also known as dynamic metamodels \cite{volpi2015development}. 
A dynamic metamodel is able to improve its fitting capability by adaptive sampling or active learning. When an adaptive sampling approach is used, the design of experiments (DoE) used for metamodel training is not defined a priori but dynamically updated, exploiting the information that becomes available during the process. The purpose of performing an adaptive DoE is to add training points anywhere it is most useful, so as to use a relatively low number of function evaluations to represent the desired function or to identify the global optimum. An example of adaptive sampling approach is the expected improvement presented in \cite{jones1998efficient}. 

In addition to dynamic metamodel, to further reduce the computational cost, 
multi-fidelity (MF) approximation methods have been developed with the aim of combining the accuracy of high-fidelity solvers with the computational cost of low-fidelity solvers \cite{serani2019adaptive}. Thus, MF metamodels are trained with a combination of high-fidelity (accurate and expensive) simulations and low-fidelity (less accurate and less expensive) simulations. Several metamodels have been used in the literature with MF data, such as non-intrusive polynomial chaos \cite{rumpfkeil2020multi}, GP \cite{coppede2019hydrodynamic}, co-kriging \cite{debaar2015-CF} and SRBF \cite{serani2019adaptive}. A survey on MF methods can be found in \cite{park2017remarks} and a discussion on the use of MF approaches can be found in \cite{giselle2019issues}. 
Adaptive sampling methods can be used also with MF metamodels, such as the upper confidence bound presented in \cite{kandasamy2016gaussian}.

In general, the performance of a metamodel is problem-dependent and determined by several concurrent issues, such as the presence of nonlinearities, the problem dimensionality, the oscillating or smooth behavior of the function, and the approach used for its training \cite{liu2018survey}. In addition, numerical solvers usually require iterative process to converge to the computational-output. The iterative process and the solution residual may affect the value of the desired computational-output with the addition of numerical noise. The presence of numerical noise can be a critical issue for the adaptive sampling/learning process.
The output-noise, if not taken into account, can deteriorate the model quality/efficiency ({\it e.g} the optimization algorithm may prematurely converge to local minima \cite{giunta1994noisy} or the adaptive sampling method may react to noise by adding many training points in noisy region, rather than selecting new points in unseen region \cite{wackers2020adaptive}). There are different strategies to deal with noise in a SBDO process, {\it e.g} 
Meliani et al. \cite{meliani2019multi} filter-out the noise by co-Kriging regression, and Wackers et al. \cite{wackers2020adaptive} use a MF-SRBF with least square regression and a MF-GPR to filter-out and assess the noise in the training set of each fidelity level. 

The MF-GPR used in this work has been applied, in authors' previous work, for a CFD-based optimization of a NACA 4-digit airfoil \cite{wackers2020adaptive} and the uncertainty quantification of an autonomous surface vehicle \cite{ficini2021uncertainty}. In both the applications the MF-GPR with 3 fidelity levels has shown better performance than the MF-GPR with 1 and 2 fidelity levels. In both the applications the training sets were affected by numerical noise especially in the lowest-fidelity, showing the importance of having a regressive formulation of the metamodel to consider the presence of numerical noise. 

The objective of the present work is to generalize and assess the robustness of the performance of the MF-GPR, through a statistical analysis, when the training sets are affected by random noise.

The MF-GPR is built as a low-fidelity metamodel corrected with metamodels of the errors/discrepancies between successive fidelity levels. 
The method filters out the numerical noise affecting the training set through regression, providing a noise-free prediction of the desired function with the associated uncertainty. 
The adaptive sampling method is based on the maximum prediction uncertainty and includes rules to automatically select the fidelity to sample, based also on the computational cost associated to each fidelity level. 

To assess the performance of the MF-GPR a set of benchmark problems is used.
The benchmark is composed by five analytical problems taken from literature \cite{clark2016engineering,wang2017generic,abdullah2019fitness,rumpfkeil2020multi}, with one and two variables. These functions are identified as representative of real world problems within the NATO AVT-331 task group on ``Goal-Driven, multi-fidelity and multidisciplinary analysis for military vehicle system level design'' \cite{beran2020comparison}. 
Each benchmark provides one, two, and three fidelity levels. Synthetic numerical noise is introduced as a normal-distributed random value with zero mean and user-defined variance not evenly distributed in the design variable space. The random values are numerically obtained with a random generator that produces random sequences of numbers.  
The noise introduces randomness in the evaluation of the objective function. Therefore, to assess the performance and the robustness of the method a statistical analysis to assess the performance of the MF-GPR with noisy training set is proposed and discussed.  
The statistical analysis is performed repeating 100 times the adaptive sampling process, varying each time the seed of the random generator. 

The performance of the MF-GPR metamodel are assessed on the benchmark problems by evaluating the accuracy of the highest-fidelity function approximation, the maximum prediction uncertainty, the verification of the metamodel predicted minimum, and the validation of the location of the global minimum and its value. The results obtained by using 1, 2, and 3 fidelity levels are compared and discussed.
 

\section{Multi-fidelity Gaussian Process Regression Model}

Given a training set $\mathcal{T}=\{{\bf x^{\prime}}_i,f({\bf x^{\prime}}_i)\}_{i=1}^{J}$, where $\mathbf{x}^{\prime} \in \mathbb{R}^{D}$ is the variables vector of dimension $D$ and $J$ is the training set size, normalizing the variables domain into a unit hypercube, the GP prediction $\tilde{f}({\bf x})$ with a constant mean and its variance ${\rm Var}[\tilde{f}({\bf x})]$ can be written as \cite{williams2006gaussian}
\begin{equation}\label{eq:GPpred}
    \tilde{f}(\mathbf{x}) =   \mathbb{E}[\mathbf{f}(\mathbf{x^{\prime}})]+ \mathbf{k}(\mathbf{x},\mathbf{x^{\prime}}) \mathbf{K}(\mathbf{x^{\prime}},\mathbf{x^{\prime}})^{-1}(\mathbf{f}(\mathbf{x^{\prime}})-\mathbb{E}[\mathbf{f}(\mathbf{x^{\prime}})]),
\end{equation} 
\begin{equation}
    {\rm Var}[\tilde{f}(\mathbf{x})] = \mathbf{K}(\mathbf{x},\mathbf{x}) - \mathbf{k}(\mathbf{x},\mathbf{x^{\prime}})^{\mathsf{T}} \mathbf{K}(\mathbf{x^{\prime}},\mathbf{x^{\prime}})^{-1}\mathbf{k}(\mathbf{x},\mathbf{x^{\prime}}),
\end{equation}
where $\mathbb{E}[\mathbf{f}(\mathbf{x^{\prime}})]$ is the expected value of $\{{f}(\mathbf{x^{\prime}}_i)\}_{i=1}^{J}$, $\mathbf{K}(\mathbf{x}^{\prime},\mathbf{x}^{\prime})$ is the covariance matrix with elements $K_{ij} = k(\mathbf{x}^{\prime}_i,\mathbf{x}^{\prime}_j)$, and $\mathbf{k}(\mathbf{x},\mathbf{x^{\prime}})$ is the covariance vector with elements $k_{i}=k(\mathbf{x},\mathbf{x}^{\prime}_i)$. Finally, $k({\cdot},{\cdot})$ is the covariance function defined as \cite{williams2006gaussian}
%
\begin{equation}
    k(\mathbf{x},\mathbf{x^{\prime}}) = \sigma_F^{2} e^{(-\boldsymbol{\gamma}^{\mathsf{T}} (\mathbf{x}-\mathbf{x^{\prime}})^{\circ 2})} + \sigma^{2}_n \delta(\mathbf{x},\mathbf{x^{\prime}}),
\end{equation}
with "$\circ$"  the Hadamard product, $\delta$ the Kronecker delta, and $\boldsymbol{\Lambda}=\{\sigma_n^2, \sigma_F^2, \boldsymbol{\gamma}\}$ the set of GP hyperparameters. Specifically, $\boldsymbol{\gamma} \in \mathbb{R}^{D}$ is the inverse length scale parameter, $\sigma_F^2$ is the signal variance, and $\sigma_n^2$ is noise variance of the training set \cite{williams2006gaussian}.

To estimate the variance associated to the noise in the training set (if not known a priori), $\sigma_n^2$ and $\sigma_F^2$ are evaluated beside $\mathbf{\gamma}$ by maximizing the log marginal likelihood $l$ \cite{williams2006gaussian} as follows
\begin{equation}\label{eq:MLE}
\boldsymbol{\Lambda}^{\star} = \{\sigma_{n}^{2,\star},\sigma_{F}^{2,\star},\boldsymbol{\gamma}^{\star}\} = {\underset{\sigma_n^2,\sigma_F^2,\boldsymbol{\gamma}}{\rm argmax}}[l],
\end{equation}
where
\begin{equation}
    l = \log p(f(\mathbf{x}^{\prime})|\mathbf{x}^{\prime}) = -\frac{J}{2} \log{2\pi} - \frac{1}{2} f(\mathbf{x}^{\prime})^T {\bf K}(\mathbf{x}^{\prime},\mathbf{x}^{\prime})^{-1} f(\mathbf{x}^{\prime}) - \frac{1}{2} \log|{\bf K}(\mathbf{x}^{\prime},\mathbf{x}^{\prime})|.
\end{equation}

The metamodel prediction uncertainty $U_{\tilde{f}}$ is here quantified as
\begin{equation}\label{eq:GPU}
    U_{\tilde{f}} = 4\sqrt{{\rm Var}[\tilde{f}({\bf x})]},
\end{equation}
yielding that $U_{\tilde{f}}$ also includes the variance associated to the noise in the training set.

Extending the metamodel prediction to an arbitrary number $N$ of fidelity levels, the MF approximation of $f_i({\bf x})$ is then built as follows \cite{wackers2020adaptive}. Given a training set $\mathcal{T}_i=\{{\bf x^{\prime}}_{j},f_i({\bf x^{\prime}}_{j})\}_{j=1}^{J_i}$ for $i=1,\dots,N$ (where $i=1$ indicates the highest-fidelity and $i=N$ indicates the lowest-fidelity), the MF approximation $\hat{f}_{i}(\mathbf{x})$ of $f_{i}(\mathbf{x})$ reads 
\begin{equation}
\label{eq:MLGeneral}
\hat{f}_{i}(\mathbf{x})\approx\widetilde{f}_N(\mathbf{x})+\sum_{k=i}^{{N}-1}\tilde{\varepsilon}_k(\mathbf{x}), 
\end{equation}
where $\tilde{{\varepsilon}}_k(\mathbf{x})$ is the inter-level error metamodel (bridge function) with an associate training set $\mathcal{E}_k=\{({\bf x^{\prime}},f_{k+1}({\bf x^{\prime}})-\hat{f}_{k}({\bf x^{\prime}}))\,|\,{\bf x^{\prime}} \in \mathcal{T}_{k+1} \cap \mathcal{T}_{k} \}$. It can be noted that Eq. \eqref{eq:MLGeneral} does not strictly require nested training sets. 

\begin{figure}[!b]
\centering
\includegraphics[width=0.5\textwidth]{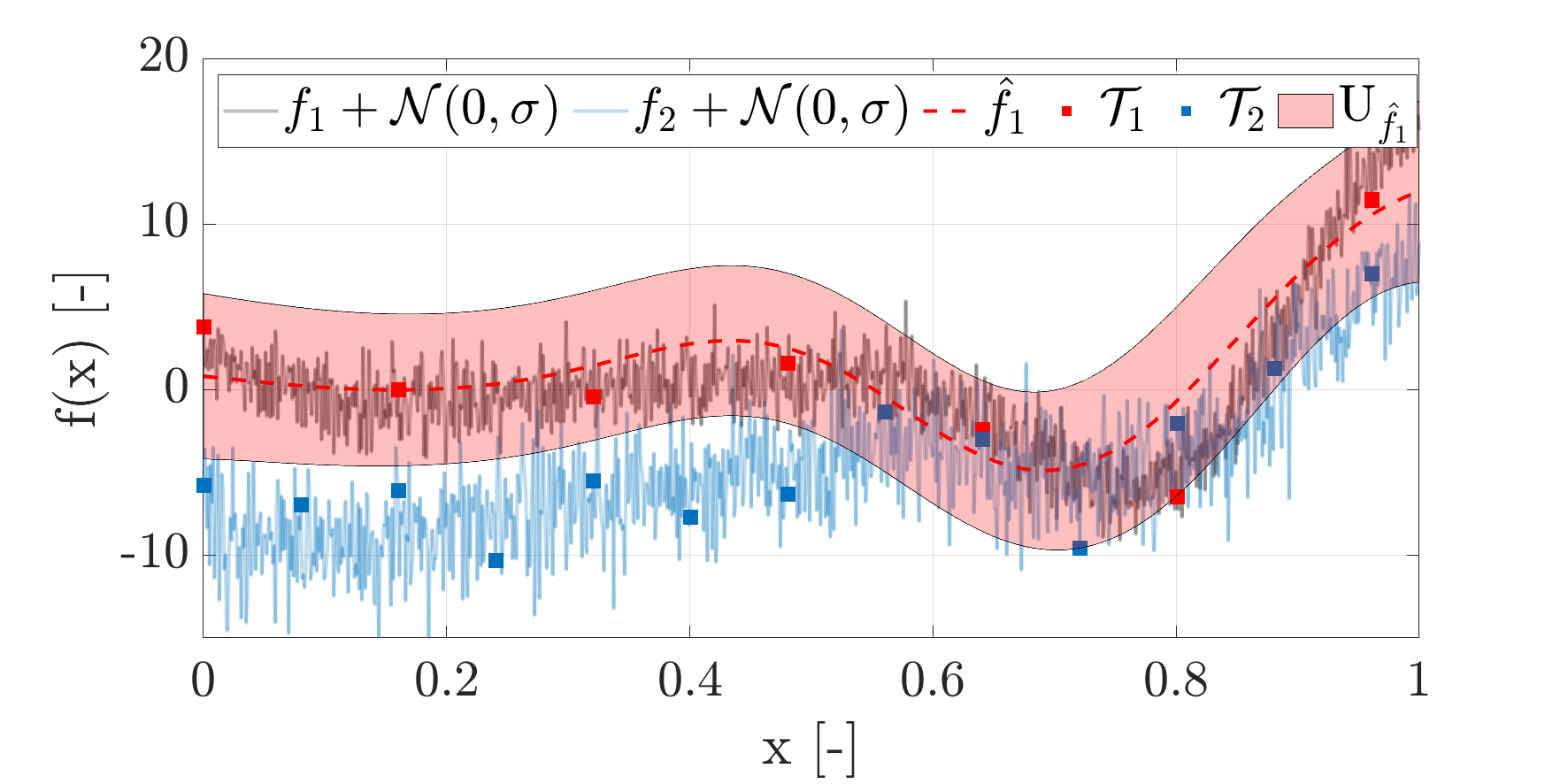}
\caption{Example of MF-GPR metamodel with two fidelity levels for noisy data.}\label{fig:MFGP_Example}
\end{figure}
Assuming the uncertainty associated to the prediction of the lowest-fidelity $U_{\widetilde{f}_N}$ and inter-level errors $U_{\tilde{\varepsilon}_k}$ as uncorrelated, the MF approximation $\hat{f}_1(\mathbf{x})$ of $f_1(\mathbf{x})$ and its uncertainty $U_{\hat{f}_1}$ read 
\begin{equation}
\label{eq:ML-MF}
\hat{f}_1(\mathbf{x})\approx\widetilde{f}_N(\mathbf{x})+\sum_{k=1}^{N-1}\tilde{\varepsilon}_k(\mathbf{x}),
~~~~~
\mathrm{and}
~~~~~
U_{\hat{f}_1}(\mathbf{x})=
\sqrt{U^2_{\widetilde{f}_N}(\mathbf{x})+\sum_{k=1}^{N-1}U^2_{{\tilde{\varepsilon}}_k}(\mathbf{x})}.
\end{equation}
\subsection{Adaptive sampling approach}
The MF-GPR metamodel is updated adding a new training point following a two steps procedure. First, the coordinates of the new training point ${\bf x}^\star$ are identified based on the maximum uncertainty prediction \cite{serani2019adaptive}, solving the single-objective maximization problem
\begin{eqnarray}\label{eq:MUAS}
{\bf x}^{\star}={\underset{{\bf x}}{\rm argmax}}[U_{\hat{f}_1}({\bf x})].
\end{eqnarray}
Second, the training set $\mathcal{T}_{i}$ is updated with the new training point $\left( \mathbf{x}^\star, f_{i}(\mathbf{x}^\star)\right)$ with $i=k,\dots,N$, where $k$ is defined as  
\begin{equation}\label{eq:choice} 
k=\mathrm{maxloc}\left[\mathbf{V}(\mathbf{x}^\star)\right], 
\qquad
\mathrm{with}
\qquad
\mathbf{V}\equiv 
\begin{Bmatrix} 
({\rm Var}[{{\tilde{\varepsilon}}_1}]- \sigma_{n, {\tilde{\varepsilon}_1}}^{2, \star} - p_1)/\beta_1\\ 
 \vdots\\ 
({\rm Var}[{{\tilde{\varepsilon}}_{N-1}}]- \sigma_{n, {\tilde{\varepsilon}_{N-1}}}^{2, \star} - p_{N-1})/\beta_{N-1} \\ 
({\rm Var}[{{\tilde{f}}_N}]- \sigma_{n, {f_N}}^{2, \star} - p_N)/\beta_N
\end{Bmatrix} 
\end{equation} 
where $\beta_{i}$ is the computational cost associated to the $i$-th fidelity level (normalized with respect to the high-fidelity one) and $p_i$ is a penalization value. The latter is used to avoid over-fitting of the training points, which would result in an ill-conditioned matrix while solving Eqs. \ref{eq:GPpred} and \ref{eq:MLE}. The penalization $p_i$ is applied only if $\mathbf{x}^\star$ lies within a minimum distance $d_{\min}=0.005$ of an already existing training point of $\mathcal{T}_i$. In such a case $p_i$ is evaluated as
\begin{equation}\label{eq:Penalization} 
p_i= \sum_{j=1}^{J_i} \frac{1}{\lVert \mathbf{x}^\star -\mathbf{x}_j^{\prime} \rVert + \tau}, 
\end{equation} 
where $\tau=0.01$ is a scalar used only to avoid null value of the denominator in Eq. \ref{eq:Penalization}. It can be noted that the subtraction of $\sigma_n^{2,\star}$ from the variance of the prediction allows to filter-out noise from the training set while selecting the fidelity level to sample. 
Wackers et al. \cite{wackers2020adaptive} have shown that the subtraction of the noise variance while selecting the fidelity level to sample allows to reduce the clusterization of the training set points, especially in the lowest fidelity. An example of MF-GPR metamodel with noisy data is shown in Fig \ref{fig:MFGP_Example}.


\section{Analytical Benchmark Problems}

\begin{figure}[!ht]
\centering
\includegraphics[width=0.45\textwidth]{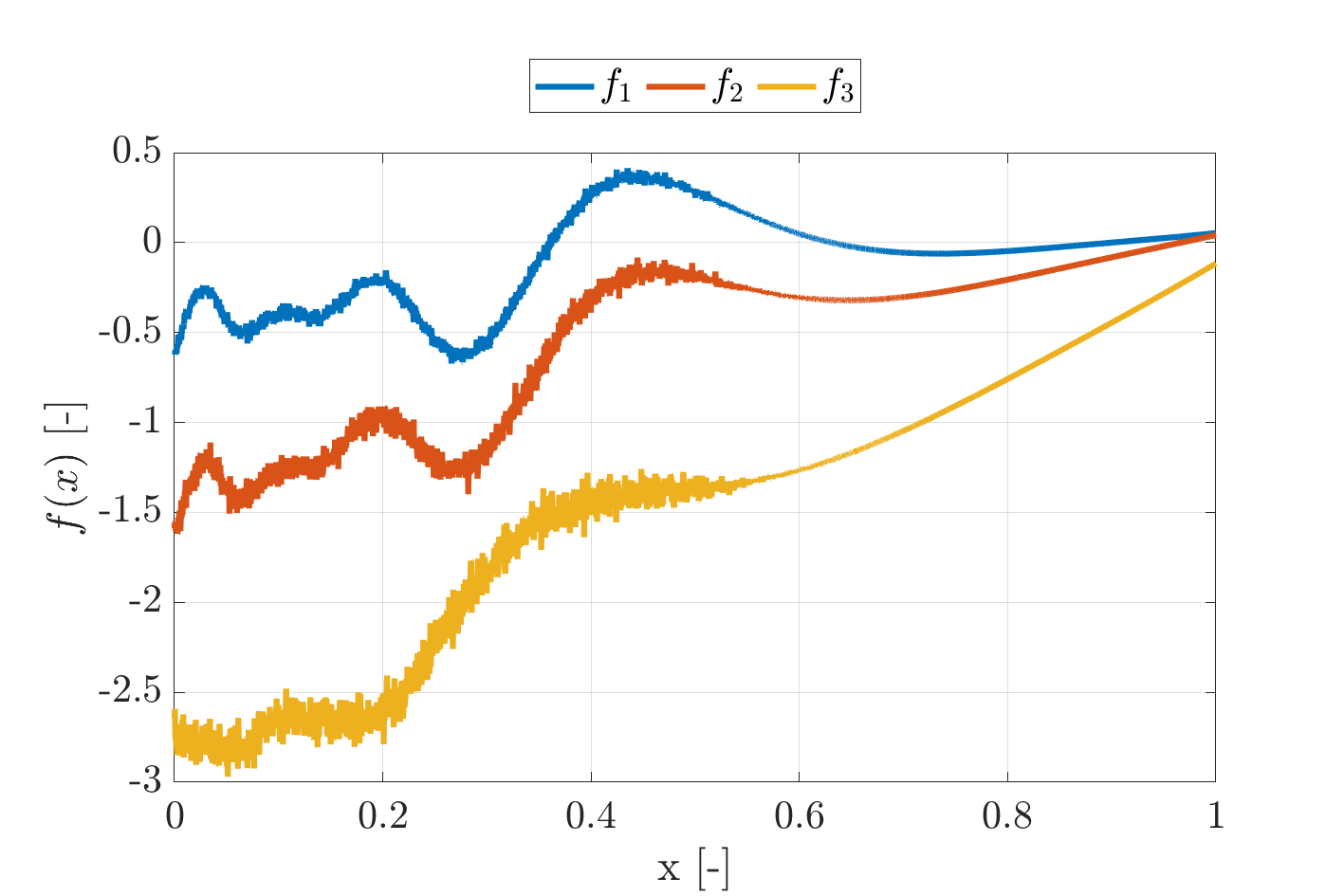}
\includegraphics[width=0.45\textwidth]{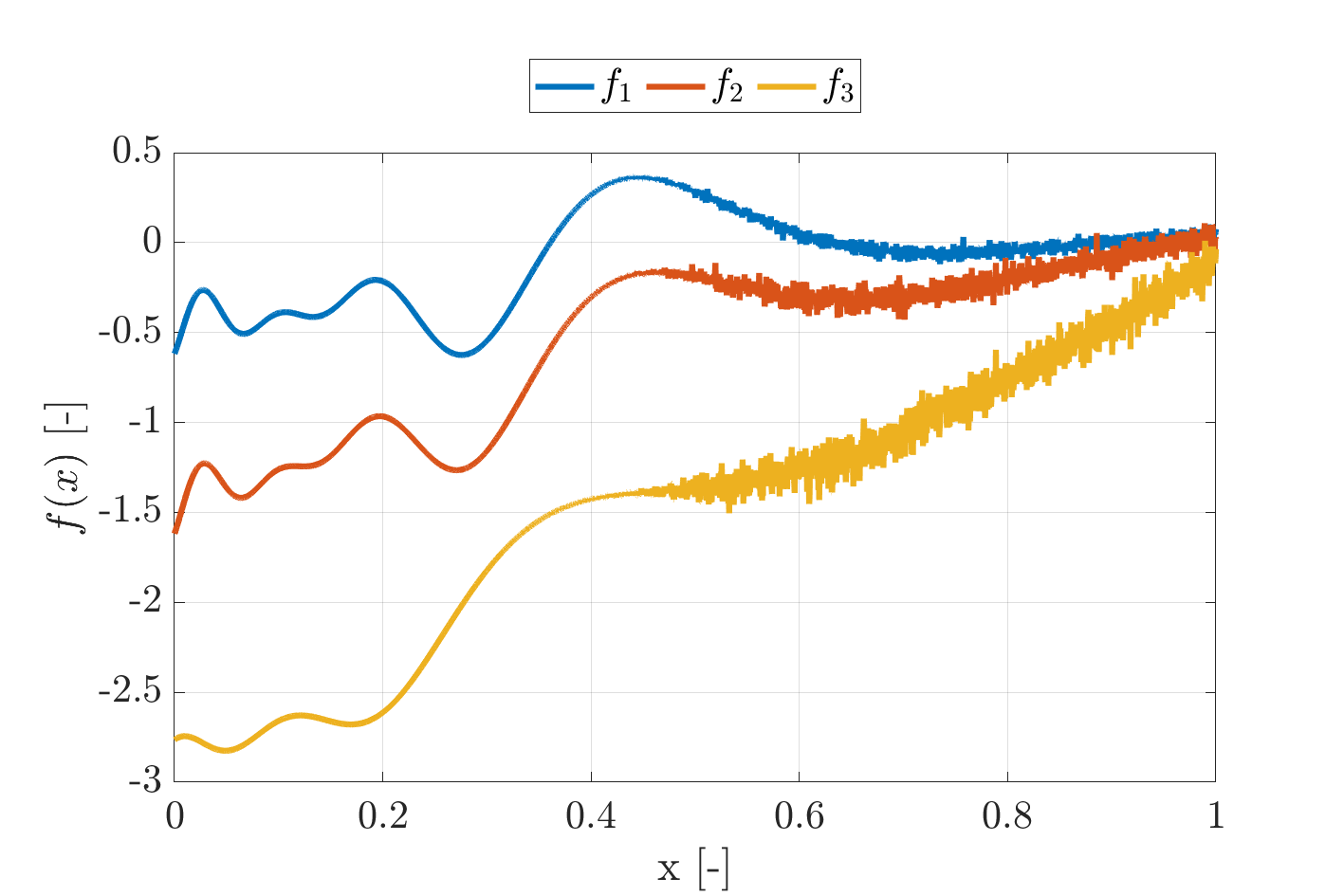}
\caption{$\mathbf{P_1}$ (left) and $\mathbf{P_2}$ (right) benchmark problems}\label{fig:Alos}
\end{figure}

\begin{figure}[!ht]
\centering
\includegraphics[width=0.32\textwidth]{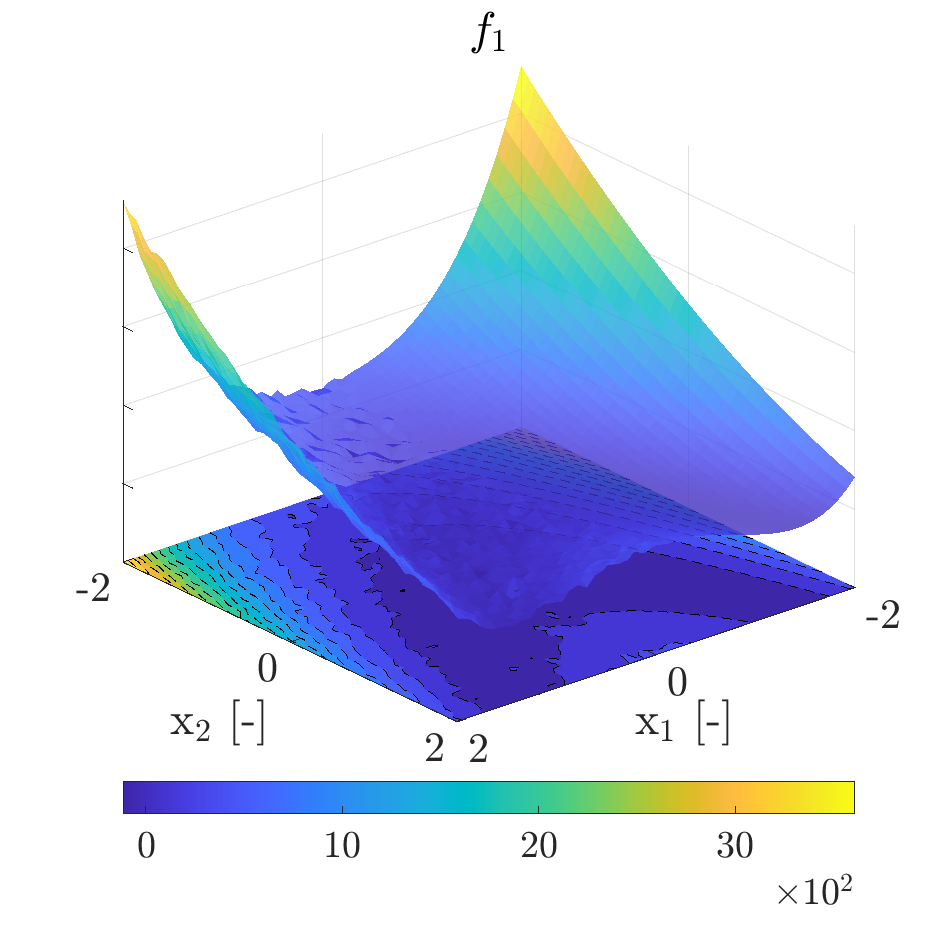}
\includegraphics[width=0.32\textwidth]{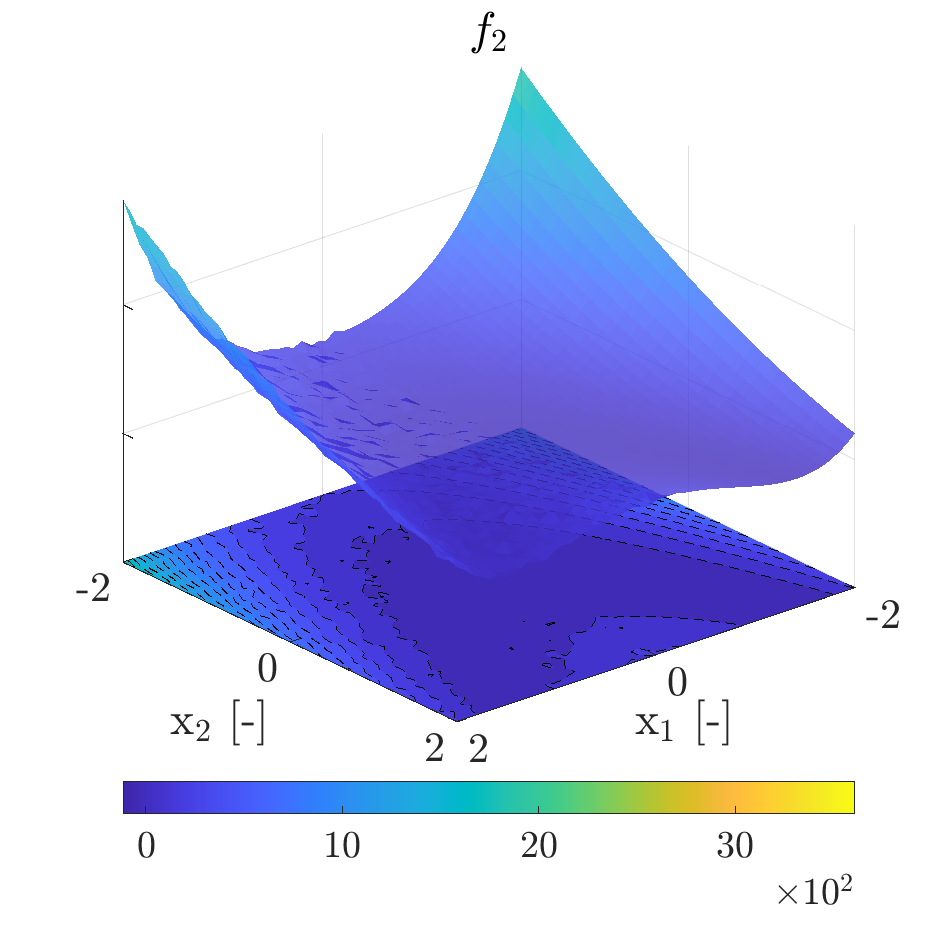}
\includegraphics[width=0.32\textwidth]{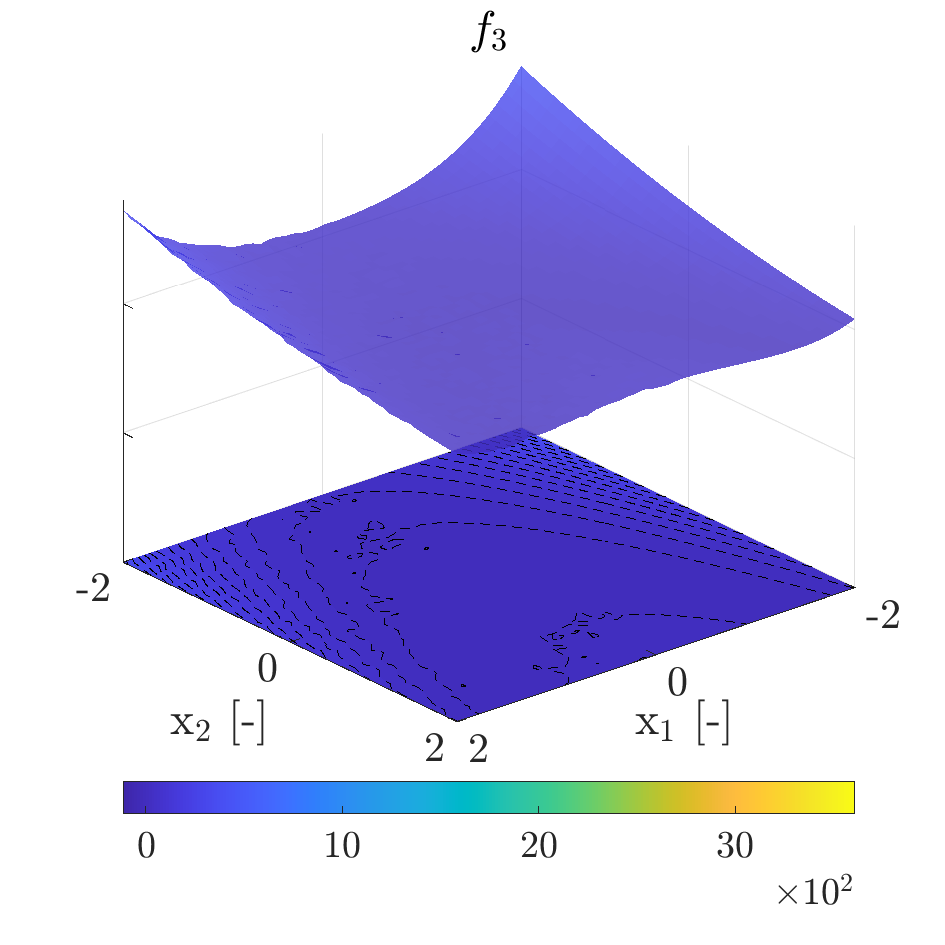}\\
\caption{$\mathbf{P_3}$ benchmark problem (from left to right $\mathbf{f_1}$, $\mathbf{f_2}$, and $\mathbf{f_3}$)}\label{fig:Rosenbrock}
\end{figure}

\begin{figure}[!ht]
\centering
\includegraphics[width=0.32\textwidth]{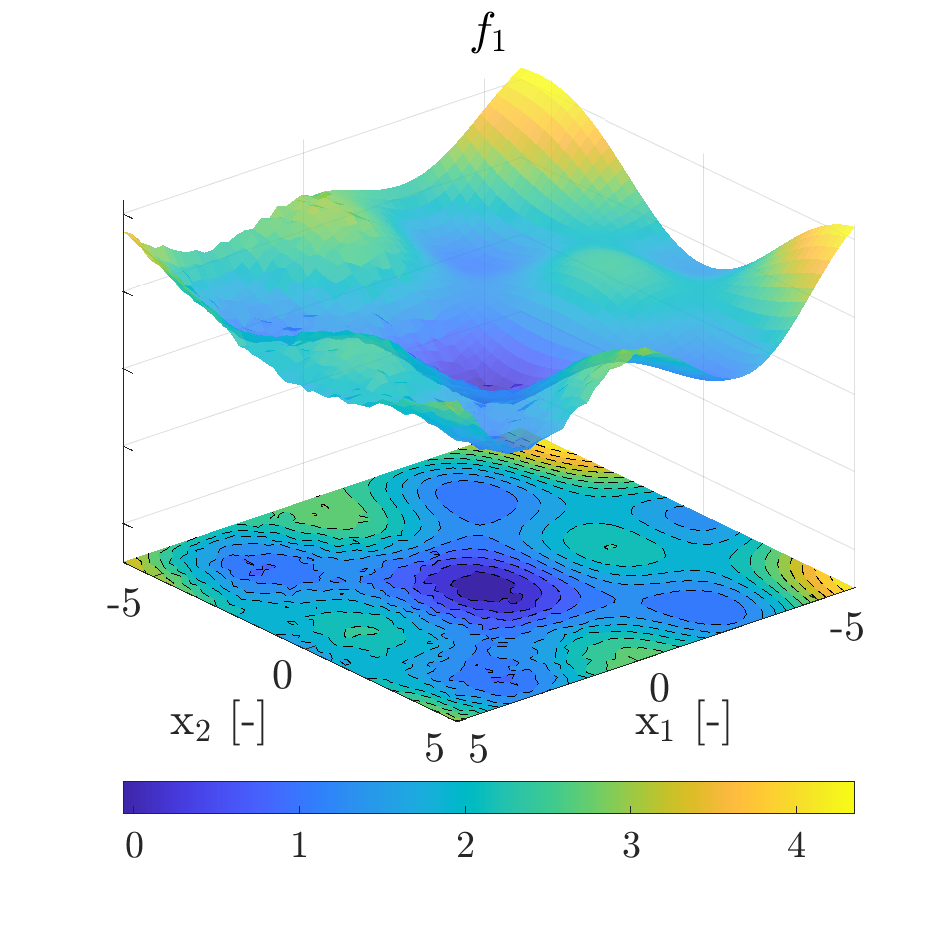}
\includegraphics[width=0.32\textwidth]{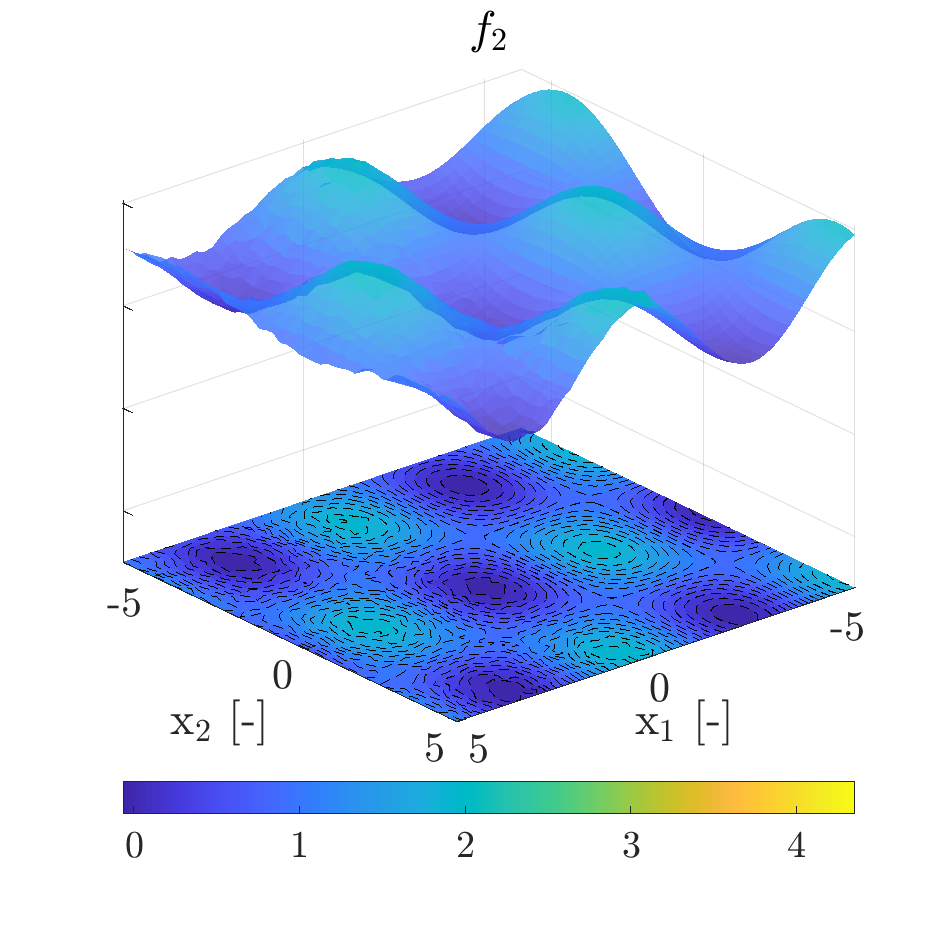}
\includegraphics[width=0.32\textwidth]{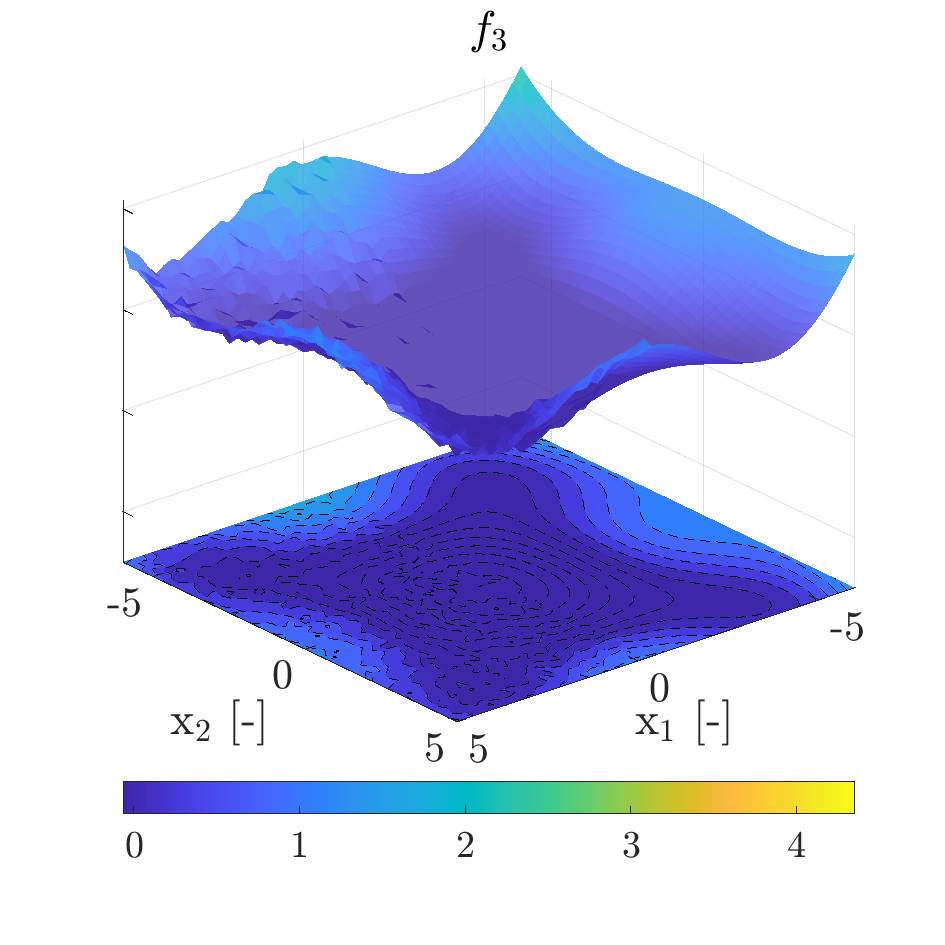}\\
\caption{$\mathbf{P_4}$ benchmark problem (from left to right $\mathbf{f_1}$, $\mathbf{f_2}$, and $\mathbf{f_3}$)}\label{fig:F11Cec}
\end{figure}

\begin{figure}[!ht]
\centering
\includegraphics[width=0.32\textwidth]{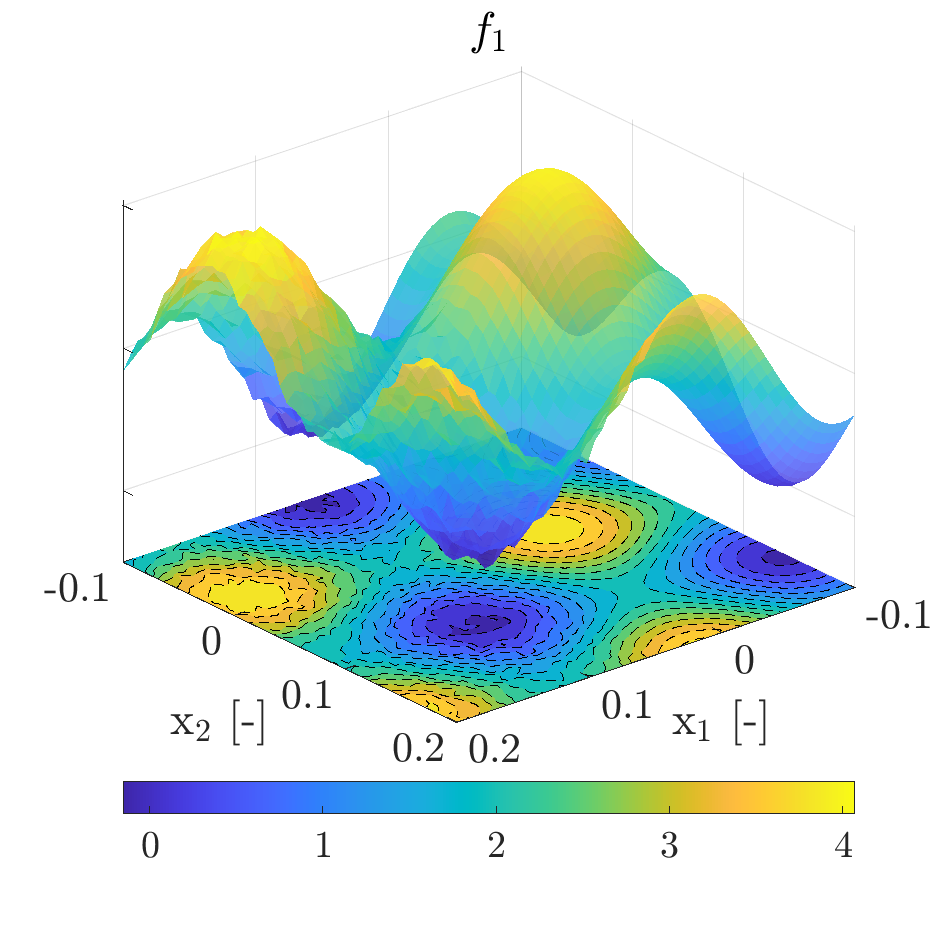}
\includegraphics[width=0.32\textwidth]{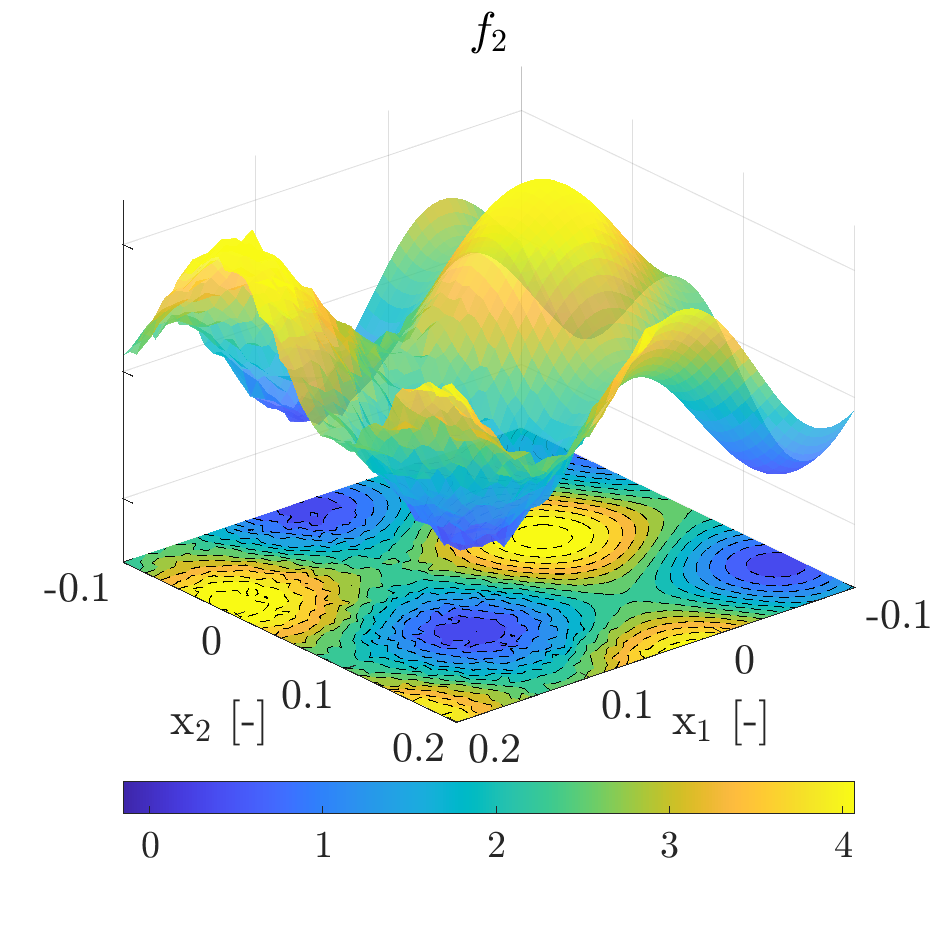}
\includegraphics[width=0.32\textwidth]{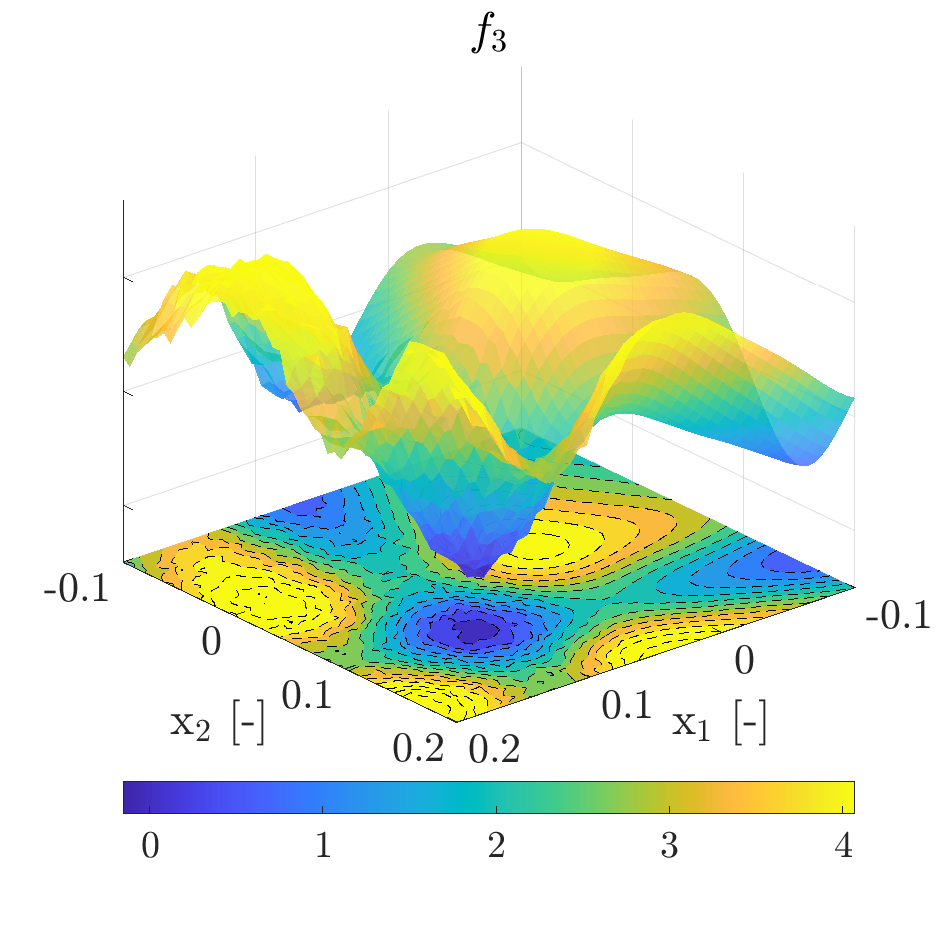}\\
\caption{$\mathbf{P_5}$ benchmark problem (from left to right $\mathbf{f_1}$, $\mathbf{f_2}$, and $\mathbf{f_3}$)}\label{fig:Rastrigin}
\end{figure}

Five analytical benchmark problems ($P$) are used to assess the MF-GPR performance and are summarized in Tab. \ref{tab:bench}. Two mono- and three bi-dimensional problems are considered. For all the benchmarks up to $N=3$ fidelities are used. For each fidelity level $i$, synthetic numerical noise is added to the analytical function and is defined as a zero mean normal distributed random value ($\eta_i\sim \mathcal{N}(0,v_i^2)$), where $v_i = a_iR_i$ with $a_i$ an arbitrary coefficient and $R_i$ the function range based on the initial training set of the $i$-th fidelity level. Furthermore, to provide a continuous but non even distribution of the noise across the design variable space, a sigmoid-like function $\lambda({\bf x})$ is used and defined as
\begin{eqnarray}\label{eq:Noise}
\lambda_1(x_1) = & 1/(1+\exp[32(x_1+0.5)])  & \text{for } P_1,P_3,P_4, \\ 
\lambda_2(x_1) = & 1/(1+\exp[-32(x_1+0.5)]) & \text{for }  P_2,  \nonumber \\
\lambda_3(x_1) = & 1/(1+\exp[-128(x_1-0.05)]) &  \text{for }  P_5. \nonumber
\end{eqnarray}
In Figs. \ref{fig:Alos}-\ref{fig:Rastrigin} are showed the benchmark problems.
The analytical function used for $P_5$ is a modified version of the one presented \cite{wang2017generic}. It is shifted to change the position of the minimum and rotated to change the properties of the function itself within the variable space as
\begin{equation}
    \mathbf{z} = R(\theta)(\mathbf{x}-\mathbf{\check{x}} ),
\end{equation}
where $\check{\mathbf{x}}=\{0.1\}_{i=1}^D$ is the reference minimum and $R(\theta)$ is the $n$-$D$ rotational matrix  \cite{aguilera2004general}, here $\theta=0.2$. 

The functions are multi-modal and challenging from the optimization viewpoint. Moreover, the presence of the noise adds a further level of complexity in terms of metamodel approximation since it can hide the real location of the minimum and the intrinsic multi-modal nature of the function.
\begin{table}[!hb]
\caption{Analytical benchmark problems}\label{tab:bench}
\centering
\begin{tabular}{crlccc}
\toprule
\textbf{Test} & \multicolumn{2}{c}{\bf Formulation} & \bf Ref. &\textbf{Domain} & \textbf{D}\\ 
\midrule
\multirow{4}{*}{$P_{1}$} & $f_1(x) =$&$ \sin(30(x-0.9)^4)\cos(2(x-0.9))+(x-0.9)/2+ \eta_1 \lambda_{1}({ x})$ & \cite{clark2016engineering} &\multirow{4}{*}{$x \in [0,1]$} &\multirow{4}{*}{1}\\
&$f_2(x) =$&$ (f_1(x)-1+x)/(1+0.25x)+ \eta_2 \lambda_{1}({ x})$& \cite{clark2016engineering} &\\
&$f_3(x) =$&$ \sin(20(x-0.87)^4)\cos(2(x-0.87))+(x-0.87)/2$& \multirow{2}{*}{--}&\\
&       &$-(2.5-(0.7x-0.14)^2)+2x +\eta_3 \lambda_{1}({ x})$ &  & \\
\midrule
\multirow{4}{*}{$P_{2}$} & $f_1(x) =$&$ \sin(30(x-0.9)^4)\cos(2(x-0.9))+(x-0.9)/2+\eta_1 \lambda_{2}({ x})$ & \cite{clark2016engineering} &\multirow{4}{*}{$x \in [0,1]$} &\multirow{4}{*}{1}\\
&$f_2(x) =$&$ (f_1(x)-1+x)/(1+0.25x)+\eta_2 \lambda_{2}({ x})$& \cite{clark2016engineering} &\\
&$f_3(x) =$&$ \sin(20(x-0.87)^4)\cos(2(x-0.87))+(x-0.87)/2$& \multirow{2}{*}{--}&\\
&       &$-(2.5-(0.7x-0.14)^2)+2x +\eta_3 \lambda_{2}({ x})$ &  & \\
\midrule
\multirow{3}{*}{$P_3$} & $f_1({\bf x}) =$ & $ \sum_{j=1}^{\mathbf{D}-1}[100 (x_{j+1} - x_j^2)^ 2 + (1 - x_j)^2] +\eta_1 \lambda_1(x_1)$ & \cite{rumpfkeil2020multi} &\multirow{3}{*}{${\bf x} \in [-2,2]$} &\multirow{3}{*}{2}\\
&   $f_2({\bf x}) =$ & $ \sum_{j=1}^{\mathbf{D}-1}[50 (x_{j+1} - x_j^2)^2 + (-2 - x_j)^2]-\sum_{j=1}^{\mathbf{D}} 0.5 x_j+\eta_2 \lambda_{1}(x_1)$& -- &
\\
&    $f_3({\bf x}) =$&$(f_1(\mathbf{x})-4-\sum_{j=1}^{\mathbf{D}}0.5x_j)/(10+\sum_{j=1}^{\mathbf{D}}0.25x_j) +\eta_3 \lambda_1(x_j)$& \cite{rumpfkeil2020multi} &\\
\midrule
\multirow{3}{*}{$P_4$} & $f_1({\bf x}) =$&$\sum_{j=1}^{\mathbf{D}}\frac{x_j^2}{25} - \prod_{j=1}^{\mathbf{D}} \cos\left(\frac{x_j}{\sqrt{j}}\right) + 1 +\eta_1 \lambda_1(x_1)$ & \cite{abdullah2019fitness} &\multirow{3}{*}{${\bf x} \in [-6,5]$}&\multirow{3}{*}{2} \\
&   $f_2({\bf x}) =$&$-\prod_{j=1}^{\mathbf{D}} \cos\left(\frac{x_j}{\sqrt{j}}\right) + 1 +\eta_2 \lambda_{1}(x_1)$& -- &
\\
&    $f_3({\bf x}) =$&$\sum_{j=1}^{\mathbf{D}}\frac{x_j^2}{20} - \prod_{j=1}^{\mathbf{D}} \cos\left(\frac{x_j}{\sqrt{j+1}}\right) -1+\eta_3 \lambda_1(x_1)$&-- &\\
\midrule
\multirow{6}{*}{$P_5$} & $f_H({\bf z}) =$ & $ \sum_{j=1}^{\mathbf{D}} (z_j^2 + 1 -\cos{(10\pi z_j)}) $ & \cite{wang2017generic} &\multirow{5}{*}{${\bf x} \in [-0.1,0.2]$} &\multirow{5}{*}{2}\\
&   $f_i({\bf z}) =$&$ f_H({\bf z}) + e_r({\bf z},\phi_i) + \eta_i \lambda_{3}(x_1)$,\,\,\,\,\,\,\,\,\,\,\,\,\,\,\,\,\,\,\,\,\,\,\ $i=1,\dots,N$ & \cite{wang2017generic} &
\\
&   $e_r({\bf z},\phi_i) =$&$\sum_{j=1}^{\mathbf{D}} a(\phi_i)\cos^2{\omega(\phi_i)z_j + b(\phi_i) + \pi} $,\,\,\,\ $i=1,\dots,N$ & \cite{wang2017generic} &
\\
&  with & $a(\phi_i)=\Theta(\phi_i)$, $\omega(\phi_i)=10\pi\Theta(\phi_i)$, $b(\phi_i)=0.5\pi\Theta(\phi_i)$ & \cite{wang2017generic} &
\\
&  and & $\Theta(\phi_i)=1-0.0001\phi_i$ & \cite{wang2017generic} &
\\
& with & $\phi = \{10000,5000,2500\}$ & - & \\
\bottomrule     
\end{tabular}
\end{table}


\section{Evaluation Metrics}
Four metrics are used to assess the performance of the MF-GPR metamodel on the benchmark problems. 
The normalized root mean squared error (NRMSE) is used to assess the accuracy of the metamodel to globally approximate the desired function and is defined as follows
\begin{equation}\label{eq:NRMSE}
{\rm NRMSE} = \frac{1}{R_1}\sqrt{\frac{1}{V}\sum_{i=1}^V[f_v({\bf x}_i) - \hat{f}({\bf x}_i)]^2},
\end{equation}
where $\{f_v(\mathbf{x}_i)\}_{i=1}^V$ is the validation set (of size $V$) evaluated by the highest-fidelity function without the noise. This metric quantify also the ability of the MF-GPR metamodel in filtering out the numerical noise. 

To quantify the accuracy of the metamodel in the predicted minimum $\hat{f}({\bf x}_{\min})$, the prediction error ($E_p$) is defined as
\begin{equation}
E_p =   \left| \dfrac{\hat{f}({\bf x}_{\min})-f_1({\bf x}_{\min})}{R_1} \right|,
\end{equation}
where $f_1({\bf x}_{\min})$ is the verified minimum by a high-fidelity evaluation without noise. 

To quantify the error in the identification the global reference minimum $f(\check{\bf x})$ in the function space, the validation error ($E_v$) is defined as
\begin{equation}
E_v = \left| \dfrac{f_1({\bf x}_{\min})-f(\check{\bf x})}{R_1} \right|,
\end{equation}
where $\check{\bf x}$ is the reference minimum.

Finally, since the benchmarks have several local minima, to quantify the effectiveness of the metamodel in identifying the position of the minimum in the variable space, the location error ($E_x$) is defined by a normalized Euclidean distance as follows
\begin{equation}
E_x=  \sqrt{\sum_{j=1}^D  \left(\frac{{x}_{\min,j}-\check{x}_j }{u_j-l_j}\right)^2},
\end{equation}
where $l_j$ and $u_j$ (for $j=1,\dots,D$) are the lower and the upper bounds of the variables domain, respectively. 

\section{Numerical Results}
The adaptive sampling procedure starts with $2D+1$ training points (for each fidelity level) located at the domain center and at the center of the domain boundaries. The metamodel based optimization as well as the solution of the minimization problem in Eqs. (\ref{eq:MLE}) and (\ref{eq:MUAS}) are based on a deterministic particle swarm optimization algorithm \cite{serani2016-ASC}. 

Since the benchmark problems are analytical functions, an artificial normalized computational cost ($cc = \sum_{i=1}^N J_i\beta_i$) is defined for each fidelity, the normalization is performed using the high-fidelity cost as reference. Thus, $\boldsymbol{\beta} = \{1\}$ is used when $N=1$, $\boldsymbol{\beta} = \{1,0.05\}$ when $N=2$, and $\boldsymbol{\beta} = \{1,0.1,0.05\}$ when $N=3$. 
When $N=1$ the metamodel is based on the highest-fidelity only, whereas when $N=2$ highest ($f_1$) and lowest-fidelities ($f_3$) are used. 
The adaptive sampling is performed considering a fixed and limited computational budget equal to $20D$. 

The coefficient sets defining the noise magnitude are set as $a_i = 0.025$ with $i=1,\dots,N$. The random values are numerically obtained with a random generator that produces random sequences of numbers. 
Since the adaptive sampling method sequentially sample the variables domain, the noise that is added in a specific point of the domain depends by the entire sampling history.
As a consequence, different methods/metamodels will see a different noise generated by a different random sequence, achieving different performance. Therefore, a statistical analysis is needed to assess the variability of the performance due to the random generator. 
Specifically, 100 runs for each benchmark problem are performed.
Finally, the NRMSE is computed on a validation set of $V=100^{D}$ points uniformly distributed in the variable's domain.
 
\begin{figure}[!t]
\centering
\subfigure[$P_1$]{\includegraphics[width=0.45\textwidth]{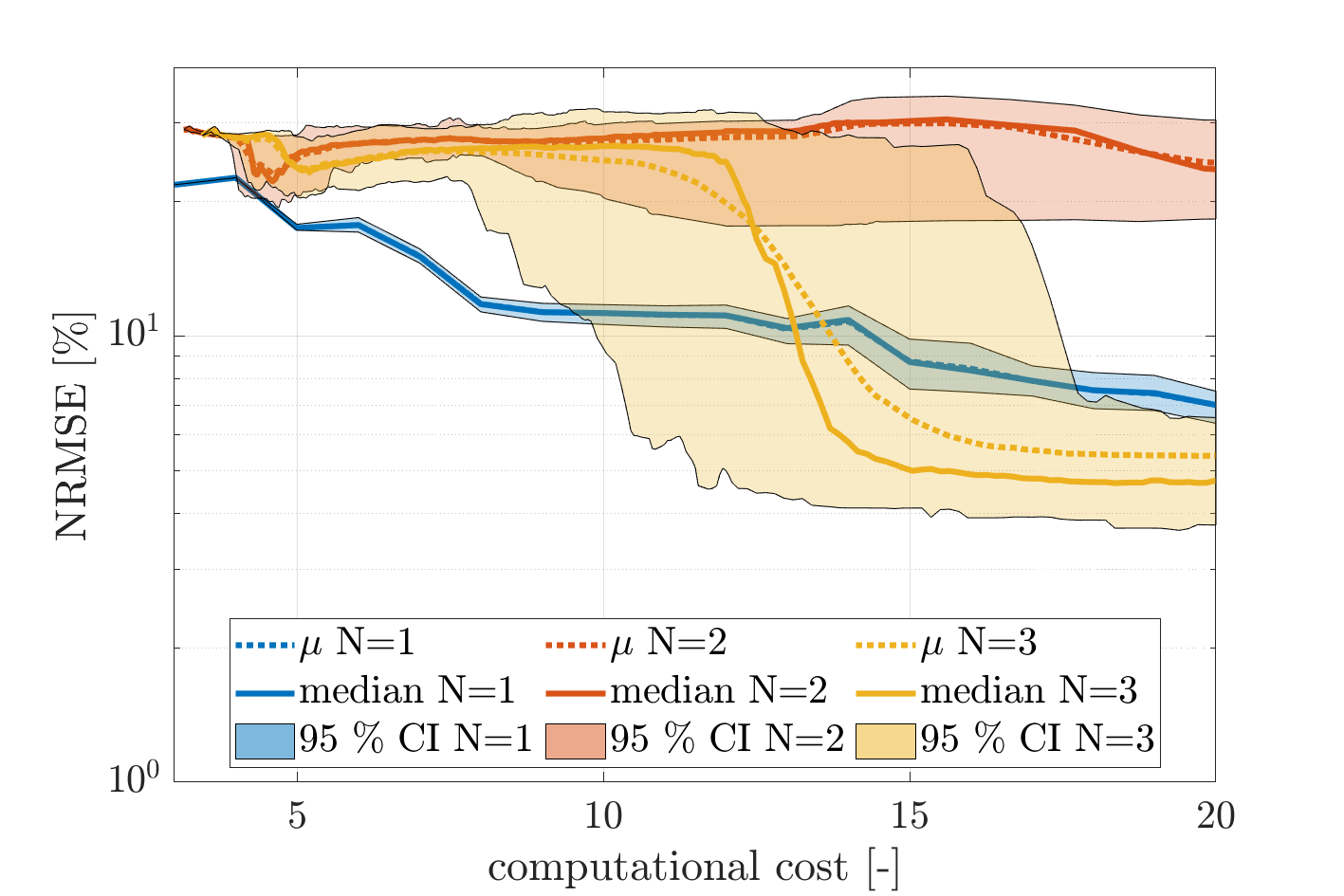}}
\subfigure[$P_2$]{\includegraphics[width=0.45\textwidth]{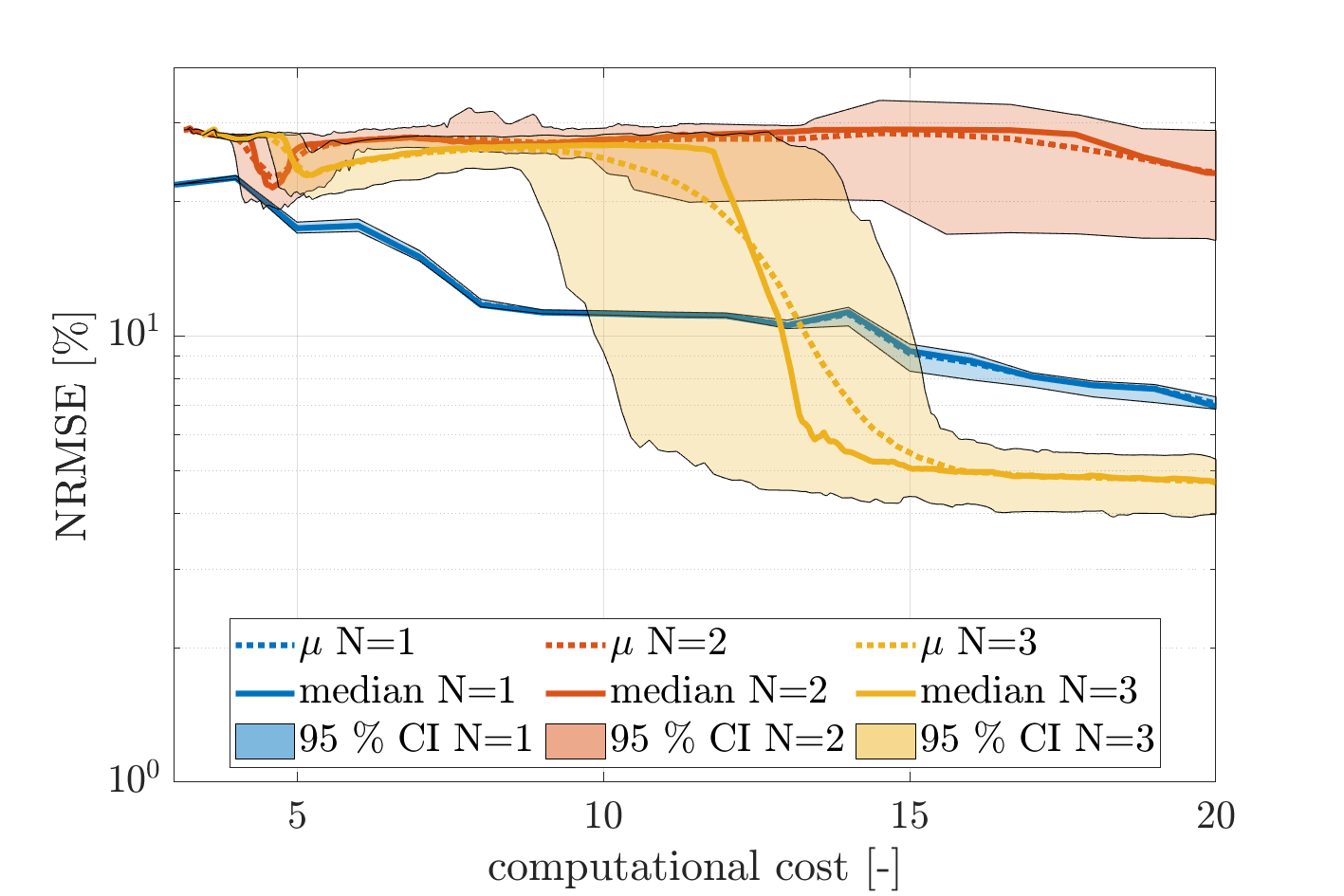}}\\
\subfigure[$P_3$]{\includegraphics[width=0.45\textwidth]{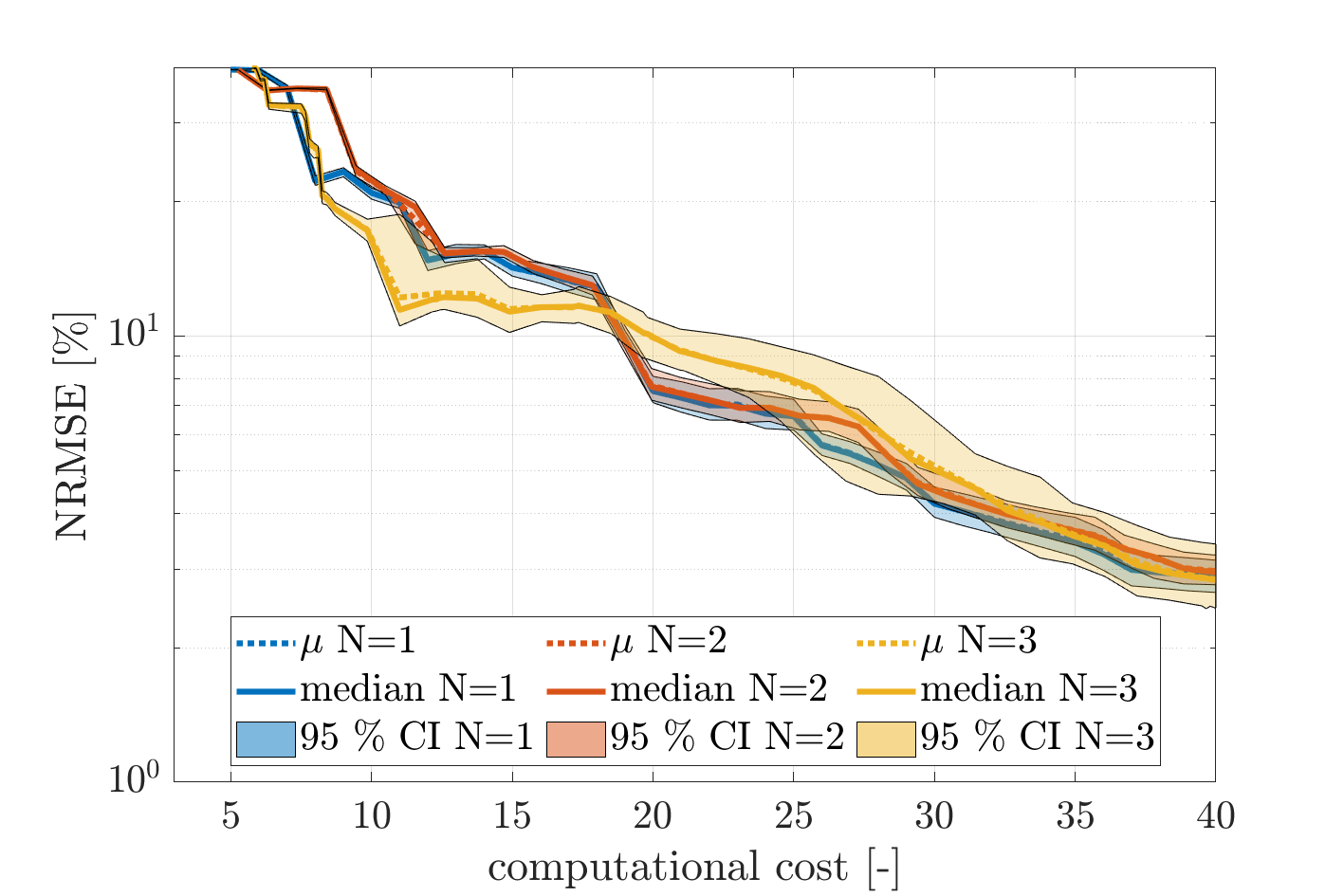}}
\subfigure[$P_4$]{\includegraphics[width=0.45\textwidth]{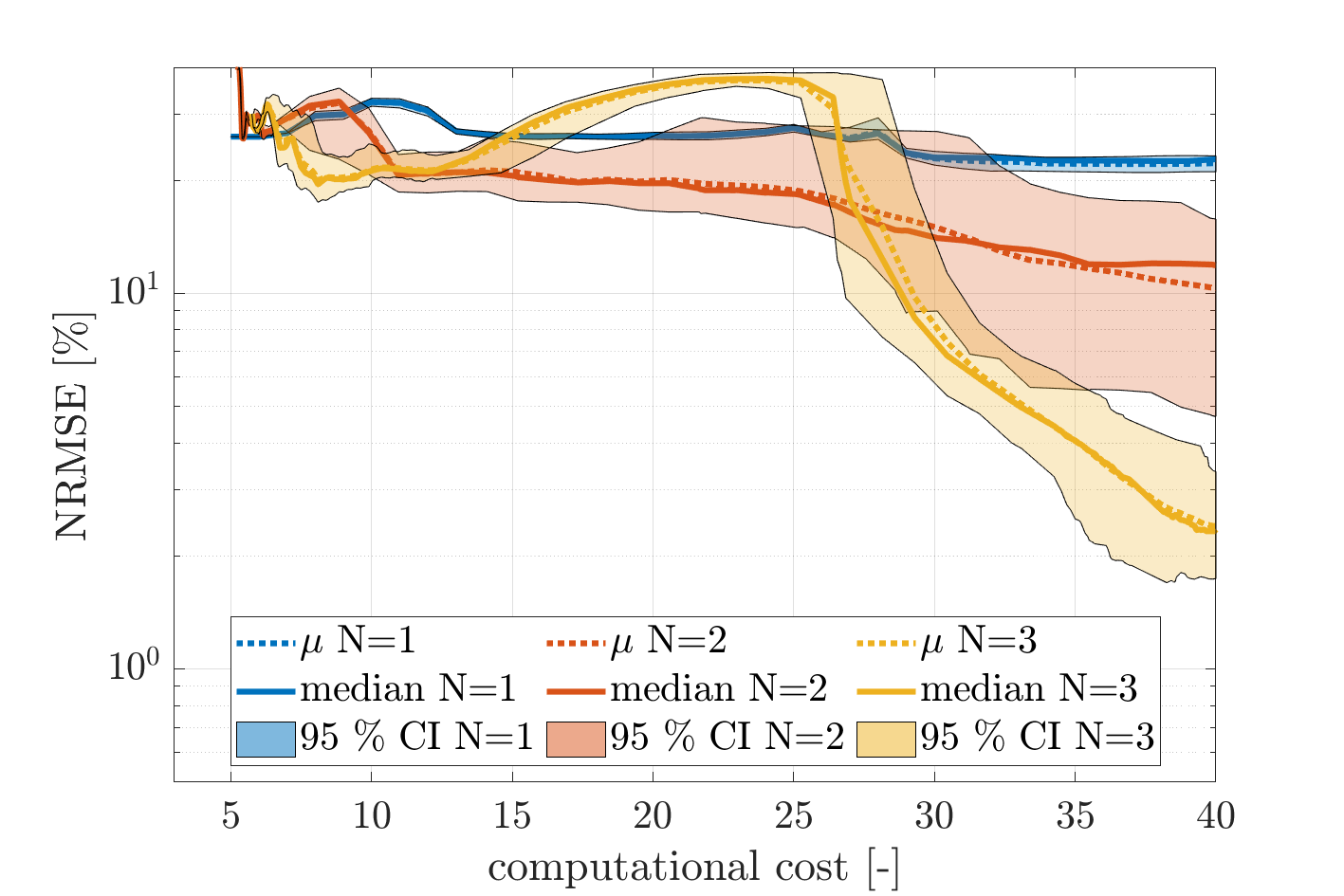}}
\subfigure[$P_5$]{\includegraphics[width=0.45\textwidth]{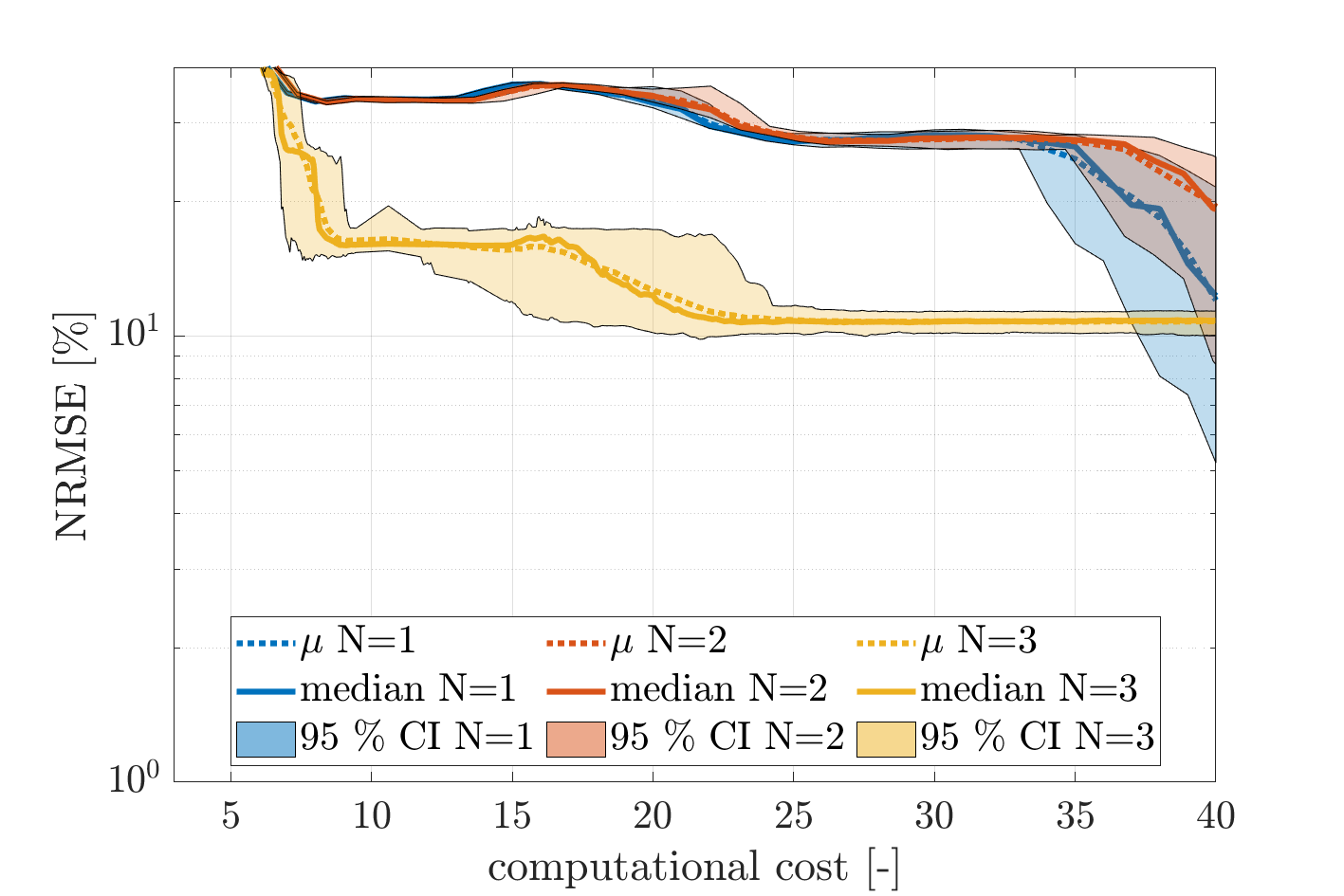}}
\caption{Convergence of the mean and the median of the NRMSE, with the 95$\%$ confidence interval, for the benchmark problems}\label{fig:NRMSE}
\end{figure}
Figure \ref{fig:NRMSE} shows the convergence of the mean and median values of the NRMSE with the 95$\%$ confidence interval, for all the benchmark problems. 
Figure \ref{fig:NRMSE}a shows that for $P_1$ the mean performance $\mu$ with $N=1$ is almost monotonically decreasing and with $N=2$ does not improve. With $N=3$ the mean performance is initially similar to $N=2$ but, with a computational cost of about 12, converges to a better value than with $N=1$. The median for $N=1$ and $N=2$ is coincident with the mean suggesting that the different convergences follow a normal distribution. For $N=3$ the median is initially higher and then lower than the mean, suggesting that the different convergences do not follow a normal distribution. Finally, the confidence interval of the convergences is smaller for $N=1$ and larger for $N=3$. It is worth noting that after a computational cost of about 18 the confidence intervals for $N=1$ and $N=3$ do not overlap, showing that the better performance of MF-GPR with $N=3$ is robust to the presence of noise. 
Figure \ref{fig:NRMSE}b shows that for $P_2$ the mean performance $\mu$ with $N=1$ is almost monotonically decreasing and with $N=2$ does not improve. With $N=3$ the mean performance is initially similar to $N=2$ but, with a computational cost of about 12, converges to a better value than with $N=1$. The median for $N=1$ and $N=2$ is coincident with the mean suggesting that the different convergences follow a normal distribution. For $N=3$ the median is initially higher and then coincident with the mean. Finally, the confidence interval of the convergences is smaller for $N=1$ and larger for $N=3$. It is worth noting that after a computational cost of about 15 the confidence interval for $N=3$ is below the one for $N=1$, showing that the better performance of MF-GPR with $N=3$ is robust to the presence of noise. 
Figure \ref{fig:NRMSE}c shows that for $P_2$ the mean performance $\mu$ with  $N=1$, 2 and 3 is monotonically decreasing. The median for $N=1$, 2 and 3 is coincident with the mean suggesting that the different convergences follow a normal distribution. Finally, the confidence interval of the convergences is larger for $N=3$. 
Figure \ref{fig:NRMSE}d shows that for $P_4$ the mean performance $\mu$ with $N=1$ does not significantly improve. With $N=2$ is monotonically decreasing. With $N=3$ the mean performance is initially worse than $N=1$ and $N=2$ but, with a computational cost of about 30 converges to a better value than with $N=1$ and $N=2$. The median for $N=1$, 2 and 3 is coincident with the mean suggesting that the different convergences follow a normal distribution. Finally, the confidence interval of the convergences is smaller for $N=1$ and larger for $N=2$. It is worth noting that after a computational cost of about 35 the confidence intervals for $N=3$ does not overlap with $N=2$. 
Figure \ref{fig:NRMSE}e shows that for $P_5$ the mean performance $\mu$ with $N=1$ and $N=2$ initially does not significantly improve, but with a computational cost of about 35 both rapidly decrease. With $N=3$ the mean performance initially is monotonically decreasing, but with a computational cost of 25 it is converged. The median for $N=1$, 2 and 3 is coincident with the mean suggesting that the different convergences follow a normal distribution. Finally, the confidence interval of the convergences is initially smaller for $N=1$ and $N=2$, but become the largest at the end of the convergences. 

Table \ref{tab:Results} summarizes the mean values of the NRMSE ($\rm \overline{NRMSE}$), of the metamodel maximum prediction uncertainty ($\overline{U}_{\hat{f}}$), of the uncertainty in the predicted minimum ($\overline{U}_{\hat{f}}({\bf x}_{\min})$), of the prediction error ($\overline{E}_p$), of the validation error ($\overline{E}_v$), of the location error ($\overline{E}_x$), and of the training set size ($\overline{J}_i$) for each fidelity level $i$. 
\begin{table}[!t]
\caption{MF-GPR metamodel summary of the average of the numerical results.}\label{tab:Results}
\centering
\begin{tabular}{ccccccccccc}
\toprule
Test & \textit{N}& ${\rm \overline{NRMSE}}$\%  & $\max\left(\overline{U}_{\hat{f}}\right)$\%$R_1$ & $\overline{U}_{\hat{f}}({\bf x}_{min})$\%$R_1$ & $\overline{E}_p\%$ & $\overline{E}_v\%$ & $\overline{E}_x \%$ & $\overline{J}_1$ & $\overline{J}_2$ & $\overline{J}_3$\\
\midrule
 \multirow{3}{*}{$P_1$} &  1 & 7.008        & 40.85          & 37.96            & 47.38                                      & 26.84            & 23.05        & 20 & - & - \\
                        &  2 & 24.46       & 48.59           & 46.18    &  \textbf{41.99}             & 36.02         & \textbf{19.62}         & 12 & - & 188\\  
                        &  3 & \textbf{5.397} & \textbf{35.38}         & \textbf{35.17}            &  67.63 &  \textbf{1.762} & 27.39 & 3 & 90 & 162\\           
\midrule
\multirow{3}{*}{$P_2$} & 1 & 7.053         & 42.26         & 39.11 &  \textbf{37.75}          & 33.49 & 25.73 &  20 & - & - \\
                        & 2 & 23.33         & 51.46         & 48.83          &  42.27          & 38.92          & \textbf{16.75}              & 12 & - & 181\\  
                        & 3 & \textbf{4.718} & \textbf{33.63} & \textbf{33.48}         &  68.52 & \textbf{0.755}          & 27.50             & 3 & 91 & 160\\
\midrule
\multirow{3}{*}{$P_3$} & 1 & 2.891 & 33.49         & 31.95            &  0.218 & 0.260 & 27.05 & 40 & - & - \\
                        & 2 & 2.972          & 31.77            & 30.31 &  0.173 & 0.209 & \textbf{17.96} & 39 & - & 40\\  
                        & 3 & \textbf{2.849}          & \textbf{26.13} & \textbf{25.36}           &  \textbf{0.130} & \textbf{0.165} & 24.85 & 33 & 52 & 53\\ 
\midrule
\multirow{3}{*}{$P_4$} & 1 & 22.23 & \textbf{70.17} & 68.22 &  7.816 & 1.755 & 2.661 & 40 & - & - \\
                       & 2 & 10.33 & 216.9 & 50.68 & 12.02 & 2.153 & 2.978 & 37 & - & 85\\  
                       & 3 & \textbf{2.401} & 189.3 & \textbf{29.69} &  \textbf{1.707} &\textbf{0.061} &\textbf{0.506} & 26 & 82 & 115\\
\midrule
\multirow{3}{*}{$P_5$}  & 1 & 12.03 & 48.63 & 26.57 & 4.392 & 1.584 & 65.37 & 40 & - & - \\
                       & 2 & 19.48 & 41.81 & 27.76 & 4.356 & 1.965 & 62.21 & 39 & - & 40\\  
                       & 3 & \textbf{10.78} & \textbf{13.10} & \textbf{12.85} & \textbf{1.717} & \textbf{1.352} & \textbf{3.208} & 10 & 131 &  349\\ 
\bottomrule   
\end{tabular}
\end{table}
MF-GPR with $N=3$ achieves the lowest ${\rm \overline{NRMSE}}$ values for all the benchmark problems considered. It may be noted that with $N=2$ the ${\rm \overline{NRMSE}}$ values are the highest for $P_1$ and $P_2$, showing that the introduction of an intermediate fidelity is beneficial in improving the MF-GPR metamodel accuracy. 
With $N=3$ the lowest value of $\max \overline{U}_{\hat{f}}$ is achieved, except for $P_4$. For $\overline{U}_{\hat{f}}(\mathbf{x}_{min})$ the lowest value is achieved with $N=3$ for all the benchmarks problems. 
With $N=2$ the lowest value of $E_p$ is achieved for $P_1$, while for $P_2$ the lowest value of $E_p$ is achieved with $N=1$. For $P_3$, $P_4$ and $P_5$ the lowest value of $\overline{E}_p$ is achieved with $N=3$.  
With $N=3$ the lowest value of $E_v$ is achieved for all the benchmarks problems. With $N=2$ the lowest value of $E_x$ is achieved for $P_1$, $P_2$ and $P_3$, while with $N=3$ the lowest value is achieved for $P_4$ and $P_5$. 
It is worth noting that $E_v$ is the lowest with $N=3$ for $P_1$ and $P_2$, while $E_p$ and $E_x$ are the highest. This can be explained because the analytical function of $P_1$ and $P_2$ present a local minimum with a value close to the value of the global minimum.

%
%

Overall, with $N=3$ there is a reduction of the number of evaluations of the high-fidelity function (up to $85\%$ for $P_1$ and $P_2$). Differently, with $N=2$ it is not guaranteed a significant reduction of the high-fidelity function evaluation.   

\section{Conclusions and Future Work}

In the context of engineering design problems, where time and computational resources are usually limited, the simulation-based design optimization has proven its ability to help designers in achieving global optimal design solutions. In the SBDO procedure are combined design modification, numerical solvers, and optimization process. When innovative design or off-design conditions are investigated, high-fidelity solvers are required. Despite the increased computational resources, the SBDO procedure with high-fidelity solvers can be very expensive from a computational viewpoint. To reduce the computational burden of the SBDO procedure adaptive single-/multi-fidelity metamodel can be used. Using few simulations, a metamodel can provide an approximate and inexpensive to evaluate model of the expensive simulations. Furthermore, the use of adaptive metamodels allow to explore the design variable space efficiently, placing new training points where is most informative. Finally, the use of multi-fidelity metamodels reduce the computational cost combining the accurate prediction of high-fidelity function evaluations with the computational cost of low-fidelity, less accurate, function evaluations. The function evaluations may be affected by numerical noise (\textit{e.g.}, due to the residuals of an iterative solver). The presence of numerical noise, if not taken into account, can deteriorate the metamodel quality/efficiency. 

In this context the assessment of the performance of an adaptive multi-fidelity metamodel based on Gaussian process regression (MF-GPR) has been presented. The assessment of the performance has been performed trough a statistical analysis on a set of five analytical benchmark problems affected by noisy function evaluation.

The MF-GPR manages an arbitrary number of fidelity levels along with random-noise affected training sets. 
The metamodel was built as the sum of a low-fidelity metamodel with metamodels of the error/discrepancy between higher-fidelity levels. The maximum prediction uncertainty was used to sequentially define new training points to adaptively refine the metamodel. The fidelity to sample was chosen based on its relative contribution to the overall prediction variance and the associated computational cost.


The benchmark problems were characterized with a synthetic numerical noise. The noise was defined as a zero mean normal distributed random value. Furthermore, to add a further level of complexity the noise was non evenly distributed in the design variable. Specifically, the noise was added in the region of the minimum, except for one mono-dimensional benchmark problem. 
Since the noise was numerically introduced using a random sequence of numbers, the noise magnitude that was added in a specific point of the domain depended by the entire history of the adaptive sampling method. As a consequence, different metamodels saw different noise generated by different random sequences.
Therefore, to assess rigorously the performance of the MF-GPR with noisy training set a statistical analysis has been done for each benchmark problems. The statistical analysis was performed repeating 100 times the adaptive sampling process of the MF-GPR varying each time the seed of the random generator.

The performances were quantified using the mean, the median and the confidence interval value obtained with the statistical analysis of four metrics, namely: the normalized root mean squared error, the prediction error, the validation error, and the location error of the minimum. 

The statistical analysis has allowed to evaluate the robustness of the MF-GPR metamodel in dealing with objective function evaluation affected by noise.
The results of the statistical analysis have shown that the MF-GPR with three fidelity levels achieved the lowest value for all the metrics considered in almost all the cases. This indicated that the MF-GPR, especially with three fidelity level, was robust in presence of noise for the accurate global approximation of the objective functions and in identifying the global minimum position and value. 
In most of the problems the MF-GPR formulation with 3 fidelity level has led to an improvement of the performance in comparison with the MF-GPR with one and two fidelities. Finally, as the number of fidelities was increased the number of high-fidelity evaluations was significantly reduced. 

Future work will focus on investigating the effects of the noise on the GP hyper-parameters evaluations. Furthermore, a larger set of benchmark problems with noisy evaluations of the objective function, considering analytical function with a larger number of variables and fidelity levels will be proposed.
Finally, the MF-GPR will be tested on a hull-form optimization used as a test case within the NATO AVT-331 task group on ``Goal-Driven, multi-fidelity and multidisciplinary analysis for military vehicle system level design'' \cite{beran2020comparison}.





\section*{Acknowledgments}
CNR-INM is grateful to Drs. Elena McCarthy and Woei-Min Lin of the Office of Naval Research for their support through the Naval International Cooperative Opportunities in Science and Technology Program. Dr. Riccardo Pellegrini is partially supported through CNR-INM project OPTIMAE. The work is conducted in collaboration with the NATO AVT-331 task group on ``Goal-driven, multi-fidelity approaches for military vehicle system-level design''.


\bibliography{biblio}

\end{document}